\RequirePackage{fix-cm}
\documentclass[smallextended]{svjour3}       
\smartqed  
\usepackage{graphicx}
\begin{document}

\title{ Study on Hilbert's Eighth Problems 
}


\author{Jinzhu Han
}


\institute{J. Han \at
              Room 1611, Wang Kezhen Building, No.52, Haidian District, Beijing 100080, People's Republic of China \\
              Tel.: +86-13314863845\\
              \email{han87654321@sina.com}           
           \and
}

\date{Received: date / Accepted: date}

\maketitle

\begin{abstract}
 In this paper, we used the principle of sieve function transformation to improve sieve method and the prime number theorem in the arithmetic sequence.For this, we proved General Riemann Hypothesis and Riemann Hypothesis to be true. further, we improved Selberg's line sieve method and proved Goldbach Conjecture and Twin Prime Conjecture to be true. However, we basically solved Hilbert's eighth problems.
\end{abstract}

\section{Introduction}
\label{intro}
  In 1900, Hilbert raised 23 famous mathematical problems,and there have been many problems that have been solved.However,Hilbert's eighth problem include Riemann Hypothesis, Goldbach conjecture, and Twin prime conjecture, that have made some progress, but it have not yet been fully resolved.  $^{[1-9]}$.

 In our paper `Transformation of sieve function', we introduced the basic principle of transformation of sieve function.$^{[11]}$.The prime number theorem in the arithmetic sequence is also improved. In this article,we will use the principle of sieve function transformation to prove General Riemann Hypothesis,and  then, we improved Selberg's line sieve method and proved Goldbach conjecture, and Twin prime conjecture. The following notations will be used.
 
  \textbf{Notations}

   $B=O(A)$,or $B\ll A$ : there be positive constant $c$ such that $|B|\leq cA$
   
   $\mu(d)$ : $ M\ddot{o}bius$ function
   
   $\nu(d)$ : the number of different prime divisors of $d$
    
   $\phi(q)$ : Euler function
   
   $\Lambda(n)$: Mangoldt function,$\Lambda(n)=\log p, n=p^{k}$, $p$ is prime
   
   $e^{\gamma}=1.78107$

\section{  The Transformation of Sieve Function }

In number theory, let there exist set  $ A=\left \{ a:a\leq x, \right \} $ be a  subset of the natural number set, which does not contain duplicate elements,  $ A_{d}=\left \{ a:a|d,a\in A \right \} $, then the  sieve function is general represented as
\begin{equation}
S(A;P(z),z)=\sum _{a\in A}\sum _{d|(a,P(z))}\mu (d)=\sum _{d|(P(z))}\mu (d)\left | A_{d} \right |
\label{eq}
\end{equation} 
where $P(z)=\prod _{p\leq z}p $ , $ |A_{d}| $ is size of set $ A_{d} $.
 In modern theory of sieve method,  $ |A_{d}| $ is represented as,
\begin{equation}
 |A_{d}| =X\frac{\omega (d)}{d}+r_{d}
\label{eq}
\end{equation} 
where $X$ is a constant, and $ \omega (d) $ is a multiplicative function, so that have a expression
\begin{equation}
S(A;P(z),z)=XW(z)+\sum _{d|(a,P(z))}\mu (d)r_{d}
\label{eq}
\end{equation} 
where
\begin{equation}
W(z)=\prod_{p|P(z)}\left(1-\frac{\omega(p)}{p}\right)
\label{eq}
\end{equation}
Let there be an anther similar sieve function
\begin{equation}
S(B;P(z),z)=\sum _{a\in B}\sum _{d|(b,P(z))}\mu (d)=\sum _{d|(P(z))}\mu (d)\left | B_{d} \right |
\label{eq}
\end{equation}
such that $|A|=|B|$,then according to the identity of the sieve function transformation,there must be transformation $S(A;P(z),z)\to  S(B;P(z),z)$,namely, $S(A;P(z),z)\sim  S(B;P(z),z)$. In our paper `Transformation of sieve function', we proved the following theorem.$^{[11]}$

\begin{theorem}
Let $S(A; P(z),z)\sim S(B; P(z),z)$, then it have
\begin{equation}
S(A; P(z),z)-S(B; P(z),z)\ll\sum _{p\leq z}\left||A_{p}|-|B_{p}|\right|+R_{0}
\label{eq}
\end{equation}
where $R_{0}$  be the error term caused by the power of primes $p^{\alpha},\alpha\geq 2$.
\end{theorem}
In the general cases,
\begin{equation}
R_{0}\ll\sum _{p\leq z}\log|A|
\label{eq}
\end{equation}
 
By this theorem, we improve the prime number theorem in the arithmetic sequence, and proved the following theorem.$^{[11]}$ 
\begin{theorem}
Let $(q,l)=1$, $l<q$ , then it have
\begin{equation}
\pi(x;q,l)=\frac{\pi(x)}{\phi(q)}+O\left(\sqrt{x}\right)
\label{eq}
\end{equation}
\end{theorem}
It is well known the function 
\begin{equation}
\psi(x;q,l)=\sum_{^{n\equiv l mod(q)}_{p\leq x}}\Lambda(n)
\label{eq}
\end{equation}
that is connected to $\pi(x;q,l)$,
\begin{equation}
\psi(x;q,l)-\pi(x;q,l)\log x\ll \sqrt{x}\log x
\label{eq}
\end{equation}
Theorem 2 equivalent following theorem.
\begin{theorem}
Let $(q,l)=1$, $l<q$ , then it have
\begin{equation}
\psi(x;q,l)=\frac{x}{\phi(q)}+O\left(\sqrt{x}\log x\right)
\label{eq}
\end{equation}
\end{theorem}
These theorem is equivalent proposition with General Riemann Hypothesis. However, let General Riemann Hypothesis be true, then Riemann Hypothesis be also true.

The principle of sieve function can be also used to improve  Selberg's sieve method.

\section{ Improving of Selberg's Sieve Method }

Selberg's sieve method have be used to research Goldbach Conjecture many years.The following theorem be proved in  Heini Halberstam and Hans-Egon Richert's book `Sieve Method',which apply to Selberg's line sieve method. $^{[10]}$
\begin{theorem}
In the line conditions
\begin{equation}
0\leq \frac{\omega(p)}{p}\leq C_{1}
\label{eq}
\end{equation}
\begin{equation}
-L\leq \sum_{w\leq p\leq z}\frac{\omega(p)\log p}{p}-\log\frac{z}{w}\leq C_{2}
\label{eq}
\end{equation}
\begin{equation}
\sum_{d<X^{\alpha}\log^{-C}X}\mu^{2}(d)3^{\nu(d)}|r_{d}|\ll \frac{X}{\log^{2}X}
\label{eq}
\end{equation}
it have
\begin{equation}
S(A;P(z),z)\leq XW(z)\left[F\left(\frac{\alpha\log X}{\log z)}\right)+O\left(\frac{L}{(\log X)^{1/14}}\right)\right]
\label{eq}
\end{equation}
and
\begin{equation}
S(A;P(z),z)\geq XW(z)\left[f\left(\frac{\alpha\log X}{\log z)}\right)+O\left(\frac{L}{(\log X)^{1/14}}\right)\right]
\label{eq}
\end{equation}
where  $L$,$C$, $C_{1}$, $C_{2}$ be calculable constants,
\begin{equation}
F(u)=\frac{2e^{\gamma}}{u},1\leq u\leq3
\label{eq}
\end{equation}
\begin{equation}
f(u)=2e^{\gamma}\frac{\log(u-1)}{u},2\leq u\leq4
\label{eq}
\end{equation}
\end{theorem}
Using principle of sieve function transformation,we improved Selberg's line sieve method, and proved following theorem.  
\begin{theorem}
In the line conditions
\begin{equation}
0\leq \frac{\omega(p)}{p}\leq C_{1}
\label{eq}
\end{equation}
\begin{equation}
-L\leq \sum_{w\leq p\leq z}\frac{\omega(p)\log p}{p}-\log\frac{z}{w}\leq C_{2}
\label{eq}
\end{equation}
it have
\begin{equation}
S(A;P(z),z)\leq XW(z)\left[F\left(\frac{\log X}{\log z)}\right)+O\left(\frac{L}{(\log X)^{1/14}}\right)\right]+R_{1}
\label{eq}
\end{equation}
and
\begin{equation}
S(A;P(z),z)\geq XW(z)\left[f\left(\frac{\log X}{\log z)}\right)+O\left(\frac{L}{(\log X)^{1/14}}\right)\right]+R_{1}
\label{eq}
\end{equation}
where 
\begin{equation}
R_{1}\ll\sum _{p\leq z}|r_{p}|+\sum _{p\leq z}\log X
\end{equation}
\end{theorem}
\begin{proof}

In above line conditions, we can set up a similar sieve function $S(A'; P(z),z)$, such that $|A|=|A'|$,and its elements be evenly distributed in set $A'$,namely,the error terms be smaller,that have
\begin{equation}
r'_{d}=|A'_{d}| -X\frac{\omega (d)}{d}\ll \log X
\label{eq}
\end{equation}
then it have
\begin{equation}
\sum_{d<X\log^{-7}X}\mu^{2}(d)3^{\nu(d)}|r'_{d}|\ll \frac{X}{\log^{2}X}
\label{eq}
\end{equation}
according to Theorem 4, we have
\begin{equation}
S(A';P(z),z)\leq XW(z)\left[F\left(\frac{\log X}{\log z)}\right)+O\left(\frac{L}{(\log X)^{1/14}}\right)\right]
\label{eq}
\end{equation}
and
\begin{equation}
S(A';P(z),z)\geq XW(z)\left[f\left(\frac{\log X}{\log z)}\right)+O\left(\frac{L}{(\log X)^{1/14}}\right)\right]
\label{eq}
\end{equation}
Because $|A|=|A'|$, $S(A; P(z),z)\sim S(A'; P(z),z)$, by Theorem 1, we have
\begin{equation}
S(A; P(z),z)-S(A'; P(z),z)\ll \sum _{p\leq z}\left||A_{p}|-|A'_{p}|\right|+R_{0}
\label{eq}
\end{equation}
\begin{equation}
|A_{p}|-|A'_{p}|=|A_{p}|-X\frac{\omega (p)}{p}+O( \log X)=r_{p}+O( \log X)
\label{eq}
\end{equation}
\begin{equation}
\sum _{p\leq z}\left||A_{p}|-|A'_{p}|\right|\ll \sum _{p\leq z}|r_{p}|+\sum _{p\leq z}\log X\ll R_{1}
\label{eq}
\end{equation}
to sum up, we have
\begin{equation}
S(A;P(z),z)\leq XW(z)\left[F\left(\frac{\alpha\log X}{\log z)}\right)+O\left(\frac{L}{(\log X)^{1/14}}\right)\right]+R_{0}+R_{1}
\label{eq}
\end{equation}
and
\begin{equation}
S(A;P(z),z)\geq XW(z)\left[f\left(\frac{\alpha\log X}{\log z)}\right)+O\left(\frac{L}{(\log X)^{1/14}}\right)\right]+R_{0}+R_{1}
\label{eq}
\end{equation}
Because
\begin{equation}
R_{0}\ll\sum _{p\leq z}\log X\ll\sum _{p\leq z}|r_{p}|+\sum _{p\leq z}\log X\ll R_{1}
\label{eq}
\end{equation}
so we have
\begin{equation}
S(A;P(z),z)\leq XW(z)\left[F\left(\frac{\alpha\log X}{\log z)}\right)+O\left(\frac{L}{(\log X)^{1/14}}\right)\right]+R_{1}
\label{eq}
\end{equation}
and
\begin{equation}
S(A;P(z),z)\geq XW(z)\left[f\left(\frac{\alpha\log X}{\log z)}\right)+O\left(\frac{L}{(\log X)^{1/14}}\right)\right]+R_{1}
\label{eq}
\end{equation}
Theorem 5 be proved.

\end{proof}
\section{Proof of Goldbach Conjecture}

Goldbach Conjecture:all evens $N>4$ can be represented as the sum of two primes.$D(N)$ be defined as the number of $N=P_{1}+P_{2}$,
\begin{equation}
D(N)=\sum _{N=p_{1}+p_{2}}1
\label{eq}
\end{equation}
We can set up a set $ A=\left \{a:a=N-p, p\leq N\right \} $,for this have
\begin{equation}
|A|=\pi(N)
\label{eq}
\end{equation}
\begin{equation}
|A_{d}|=\pi(N;d,N)
\label{eq}
\end{equation}
then it have
\begin{equation}
S(A;P(z),\sqrt{N})=\sum _{d|(P(\sqrt{N}))}\mu (d)\pi(N;d,N)
\label{eq}
\end{equation}
and
\begin{equation}
D(N)= S(A;P(z),\sqrt{N})+O(\sqrt{N})
\label{eq}
\end{equation}

let us prove $D(N)>0$ for $N$ be enough large even, then Goldbach Conjecture be true.However,we proved following theorem.
\begin{theorem}
for $N$ be enough large even,we have
\begin{equation}
S(A;P(z),\sqrt{N})\geq 0.5C(N)\frac{N}{\log^{2} N}
\label{eq}
\end{equation}
where
\begin{equation}
C(N)=\prod _{2<p\leq z}\left(1-\frac{1}{(p-1)^{2}}\right)\prod _{2< p|N}\frac{p-1}{p-2}>\frac{1}{2}
\label{eq}
\end{equation}
\end{theorem}
\begin{proof}

According to basic relationship,
\begin{equation}
S(A;P(z),\sqrt{N})=S(A;P(z),\sqrt[3]{N})-\sum_{\sqrt[3]{N}\leq p<\sqrt{N}}S(A_{p};P(z),p)
\label{eq}
\end{equation} 
then we have
\begin{equation}
S(A;P(z),\sqrt{N})\geq S(A;P(z),\sqrt[3]{N})-\sum_{\sqrt[3]{N}\leq p<\sqrt{N}}S(A_{p};P(z),\sqrt[3]{N})
\label{eq}
\end{equation}
where
\begin{equation}
|A|=\pi(N)\sim \frac{N}{\log N}
\label{eq}
\end{equation}
\begin{equation}
|A_{d}|=\pi(N;d,N)=\frac{N}{\phi(d)\log N}+r_{d}
\label{eq}
\end{equation}
\begin{equation}
|A_{pd}|=\pi(N;pd,N)=\frac{N}{\phi(pd)\log N}+r_{pd}=\frac{N}{(p-1)\phi(d)\log N}+r_{pd}
\label{eq}
\end{equation}
by Theorem 2, for $(d,N)=1$, even if $N\geq d$,we still have
\begin{equation}
r_{d}=\pi(N;d,N)-\frac{N}{\phi(d)\log N}\ll\sqrt{N}
\label{eq}
\end{equation}
\begin{equation}
r_{pd}=\pi(N;pd,N)-\frac{N}{\phi(pd)\log N}\ll\sqrt{N}
\label{eq}
\end{equation}
further, according to Theorem 5, we have
\begin{equation}
S(A;P(z),\sqrt[3]{N})\geq 2C(N)\frac{N}{\log^{2} N}\left[f\left(\frac{\log N}{\log \sqrt[3]{N})}\right)+O\left(\frac{L}{(\log N)^{1/14}}\right)\right]+R_{2}
\label{eq}
\end{equation}
and
\begin{equation}
S(A_{p};P(z),\sqrt[3]{N})\leq 2C(N)\frac{N}{(p-1)\log^{2} N}\left[F\left(\frac{\log N}{\log \sqrt[3]{N})}\right)+O\left(\frac{L}{(\log N)^{1/14}}\right)\right]+R_{3}
\label{eq}
\end{equation}
where 
\begin{equation}
R_{2}\ll \sum _{^{p\leq \sqrt[3]{N}}_{(p,N)=1}}|r_{p}|+\sum _{p\leq \sqrt[3]{N}}\log N\ll N^{5/6}\log N
\end{equation}
\begin{equation}
R_{3}\ll \sum _{^{p_{i}\leq \sqrt[3]{N}}_{(pp_{i},N)=1}}|r_{p_{i}p}|+\sum _{p\leq \sqrt[3]{N}}\log N\ll N^{5/6}\log N
\end{equation}
for $N$ be enough large even, $N^{5/6}\log N$ be far less than $N/\log^{2}N$,
\begin{equation}
f(3)=\frac{2}{3}e^{\gamma}\log 2
\end{equation}
\begin{equation}
F(3)=\frac{2}{3}e^{\gamma}
\end{equation}
\begin{equation}
\frac{L}{(\log N)^{1/14}}\to 0
\end{equation}
so we have
\begin{equation}
S(A;P(z),\sqrt[3]{N})\geq \frac{4}{3}e^{\gamma}\log 2C(N)\frac{N}{\log^{2} N}
\label{eq}
\end{equation}
\begin{equation}
S(A_{p};P(z),\sqrt[3]{N})\leq \frac{4}{3}e^{\gamma}C(N)\frac{N}{(p-1)\log^{2} N}
\label{eq}
\end{equation}
and
\begin{equation}
S(A;P(z),\sqrt{N})\geq \frac{4}{3}e^{\gamma}C(N)\frac{N}{\log^{2} N}\left[\log 2-\sum_{\sqrt[3]{N}\leq p\leq \sqrt{N}}\frac{1}{p-1}\right]
\label{eq}
\end{equation}
because for$N$ be enough large even
\begin{equation}
\sum_{\sqrt[3]{N}\leq p\leq \sqrt{N}}\frac{1}{p}\leq 1.09(\log 3-\log 2)
\label{eq}
\end{equation}
\begin{equation}
\sum_{\sqrt[3]{N}\leq p\leq \sqrt{N}}\frac{1}{p-1}\leq 1.1(\log 3-\log 2)
\label{eq}
\end{equation}
so we have
\begin{equation}
S(A;P(z),\sqrt{N})\geq \frac{4}{3}(2.1\log 2-1.1\log 3)e^{\gamma}C(N)\frac{N}{\log^{2} N}\geq 0.5C(N)\frac{N}{\log^{2} N}
\label{eq}
\end{equation}
Theorem 6 be proved.

For this,for $N$ be enough large even, $D(N)>0$. Goldbach Conjecture be true.
\end{proof}

\section{Proof of Twin Primes  Conjecture}

The Twin Primes Conjecture:the number of twin primes be infinite.

In the number theory, $T(N)$ be defined as the number of $P_{1}-P_{2}=2$,
\begin{equation}
T(N)=\sum_{^{P_{1}-P_{2}=2}_{p_{1},p_{2}\leq N}}1
\label{eq}
\end{equation}
We can  introduce the following  sieve function.
\begin{equation}
S(B;P(z),\sqrt{N})=\sum _{d|(P(\sqrt{N}))}\mu (d)\pi(N;d,2)
\label{eq}
\end{equation}
then it have
\begin{equation}
T(N)= S(B;P(z),\sqrt{N})+O(\sqrt{N})
\label{eq}
\end{equation}

We can use similar above method  to prove the  following theorem.
\begin{theorem}
for $N$ be enough large even,we have
\begin{equation}
S(B;P(z),\sqrt{N})\geq 0.5\frac{N}{\log^{2} N}
\label{eq}
\end{equation}
\end{theorem}
  
\begin{proof}
Similarly,we have
\begin{equation}
S(B;P(z),\sqrt{N})\geq S(B;P(z),\sqrt[3]{N})-\sum_{\sqrt[3]{N}\leq p<\sqrt{N}}S(B_{p};P(z),\sqrt[3]{N})
\label{eq}
\end{equation}
Using similar method,for $N$ be enough large even, we can prove that 
\begin{equation}
S(B;P(z),\sqrt[3]{N})\geq \frac{4}{3}e^{\gamma}\log 2\frac{N}{\log^{2} N}
\label{eq}
\end{equation}
\begin{equation}
S(B_{p};P(z),\sqrt[3]{N})\leq \frac{4}{3}e^{\gamma}\frac{N}{(p-1)\log^{2} N}
\label{eq}
\end{equation}
and
\begin{equation}
S(B;P(z),\sqrt{N})\geq \frac{4}{3}e^{\gamma}\frac{N}{\log^{2} N}\left[\log 2-\sum_{\sqrt[3]{N}\leq p\leq \sqrt{N}}\frac{1}{p-1}\right]
\label{eq}
\end{equation}
for $N$ be enough large even,we have
\begin{equation}
S(B;P(z),\sqrt{N})\geq 0.5\frac{N}{\log^{2} N}
\label{eq}
\end{equation}
Theorem 7 be proved.
\end{proof}
According to Theorem 7, for $N\to\infty$,
\begin{equation}
T(N)=S(B;P(z),\sqrt{N})+O(\sqrt{N})\to \infty
\label{eq}
\end{equation}
therefore, the number of twin primes be infinite, so Twin Primes Conjecture be true.

\section{Conclusion}

In this paper,we use the method of sieve function transformation to improve the prime number theorem in arithmetic sequence, from this,we proved General Riemann Hypothesis is true.Further, we proved Goldbach Conjecture and Twin Primes Conjecture to be true.However,Hilbert's eighth problems be solved.

\end{document}